\documentclass[a4paper]{article}
\usepackage{latexsym}
\usepackage{amssymb}
\usepackage{amsmath}
\usepackage{url}
\usepackage[T2A,T1]{fontenc}

\newtheorem{Theorem} {Theorem} [section]
\newtheorem{Proposition} [Theorem] {Proposition}
\newtheorem{Lemma} [Theorem] {Lemma}

\newcommand{\Proof}{ \noindent{\bf Proof.}\quad }
\newcommand{\qed}{\hfill$\Box$\medskip}

\newcommand{\restrictedto}[1]{{\big|_#1}}
\newcommand{\Ff}{{\mathbb F}}
\newcommand{\eps}{\varepsilon}
\newcommand{\adj}{\sim}
\newcommand{\nadj}{\not\sim}
\newcommand{\Eps}{{\rm E}}
\newcommand{\Tau}{{\rm T}}
\newcommand{\<}{\langle}
\renewcommand{\>}{\rangle} 
\renewcommand{\phi}{\varphi} 
\def\cH{{\mathcal H}} 
\def\bD{{\bf D}} 
\def\bL{{\bf L}} 
\def\bP{{\bf P}} 
\newcommand{\downstrut}{\rule[-0.2cm]{0cm}{0cm}}
\newcommand{\maysplit}{\allowbreak}
\newcommand\Ffss[2]{\Ff_{\hspace{-1pt}#1}^{\hspace{1pt}#2}}
\newcommand\subsupr[2]{_{#1}^{\hspace{1pt}#2}}
\newcommand\mysqueeze[2]{%
\newdimen\origspacing%
\newdimen\newspacing%
\origspacing=\fontdimen2\font%
\setlength{\newspacing}{1\origspacing}%
\addtolength{\newspacing}{-#1}%
\fontdimen2\font=\newspacing%
{#2}%
\fontdimen2\font=\origspacing}

\title{Strongly regular graphs satisfying the 4-vertex condition}
\author{

A. E. Brouwer\footnotemark[1]

\& F. Ihringer\footnotemark[2]

\& W. M. Kantor\footnotemark[1]}
\date{30 Jun 2021}

\begin{document}
\maketitle
\footnotetext[1]{Retired}
\footnotetext[2]{Dept.~of Mathematics:~Analysis, Logic and Discrete
Math., Ghent University, Belgium.
E-mail: {\tt ferdinand.ihringer@ugent.be}}

\begin{abstract}
We survey the area of strongly regular graphs satisfying
the 4-vertex condition and find several new families.
We describe a switching operation on collinearity graphs of polar spaces
that produces cospectral graphs. The obtained graphs satisfy
the 4-vertex condition if the original graph belongs to
a symplectic polar space.
\end{abstract}

\section{Introduction}
In this note we look at graphs with high combinatorial regularity,
where this regularity is not an obvious consequence of properties
of their group of automorphisms.

%
A graph $\Gamma$ is said to satisfy the {\em $t$-vertex condition} if,
for all triples $(T,x_0,y_0)$ consisting of a $t$-vertex graph $T$
together with two distinct distinguished vertices $x_0,y_0$ of $T$,
and all pairs of distinct vertices $x,y$ of $\Gamma$, the number of
isomorphic copies of $T$ in $\Gamma$, where the isomorphism maps
$x_0$ to $x$ and $y_0$ to $y$, does not depend on the choice of the pair
$x,y$ but only on whether $x,y$ are adjacent or nonadjacent. 

This concept was introduced by Hestenes \& Higman \cite{HestenesHigman71}
(who refer to the unpublished Sims \cite{Sims-unpub})
in order to study rank 3 graphs.
Clearly, a rank 3 graph satisfies the $t$-vertex condition for all $t$.
If the graph $\Gamma$ satisfies the $t$-vertex condition, where $\Gamma$
has $v$ vertices and $3 \le t \le v$, then $\Gamma$ also satisfies the
$(t-1)$-vertex condition. A graph satisfies the 3-vertex condition
if and only if it is strongly regular (or complete or edgeless).
%
%
It satisfies the $v$-vertex condition if and only if it is rank 3.
Thus, we get a hierarchy of conditions of increasing 
strength between strongly regular and rank 3.

The present paper will focus almost exclusively on the case $t=4$.
A simple criterion for the 4-vertex condition is given in
Proposition \ref{prop:sims}.
Previously not many graphs were known 
that satisfy the 4-vertex condition without being rank 3.
Here we survey the known examples and give several new constructions. 
One of our constructions proceeds by switching symplectic graphs
(see Section \ref{sec:switched_col}). As a consequence we find

\begin{Theorem} \label{thm:4vtxsrgs}
 For $v \ge 4$ there are at least $\lfloor v^{1/6} \rfloor !$
 strongly regular graphs of order at most $v$ satisfying
 the $4$-vertex condition.
\end{Theorem}

It follows that among all non-isomorphic strongly regular graphs
of order at most $v$ that satisfy the $4$-vertex condition
the fraction that is determined by their spectrum goes to 0
when $v$ goes to infinity.

%

\section{The 4-vertex condition} \label{sec:def}

{\smallskip\footnotesize
A graph of order $v$ is called {\em strongly regular} with parameters
$(v, k, \lambda, \mu)$ if it is neither complete nor edgeless, 
each vertex has degree $k$, any two adjacent vertices have exactly
$\lambda$ common neighbors, and any two non-adjacent vertices have exactly
$\mu$ common neighbors.

A graph with vertex set $V$ has {\em rank $r$} if
its automorphism group is transitive on $V$ and
has exactly $r$ orbits on $V \times V$.
Rank 3 graphs are strongly regular.

If $x$ is a vertex of the graph $\Gamma$, then the {\em local graph}
$\Gamma(x)$ of $\Gamma$ at $x$ is the induced subgraph in $\Gamma$
on the neighborhood of $x$. We say that $\Gamma$ is {\em locally} P 
when all local graphs of $\Gamma$  have property P.
If $\Gamma$ is strongly regular, then its {\em 1st subconstituent}
(at a vertex $x$) is the local graph at $x$, while its {\em 2nd subconstituent}
(at $x$) is the induced subgraph on the non-neighborhood of $x$.
If $xy$ is an edge (resp. nonedge) in $\Gamma$, then the subgraph
induced on $\Gamma(x) \cap \Gamma(y)$ is called a $\lambda$-graph
(resp. $\mu$-graph).

See \cite{BrouwerVM21} for further information about
strongly regular graphs.

\par\medskip}

Details on the parameters of graphs satisfying the 4-vertex condition
are given in \cite{HestenesHigman71}. In particular, we have the following
simple criterion for the 4-vertex condition:

\begin{Proposition} {\rm (Sims \cite{Sims-unpub})} \label{prop:sims}
A strongly regular graph $\Gamma$ with parameters
$(v,k,\maysplit\lambda,\maysplit\mu)$
satisfies the 4-vertex condition, with parameters $(\alpha,\beta)$,
if and only if the number of edges in $\Gamma(x) \cap \Gamma(y)$
is $\alpha$ (resp. $\beta$) whenever the vertices $x,y$ are adjacent
(resp. nonadjacent).
In this case, $k\big (\binom{\lambda}{2} - \alpha\big) = \beta(v-k-1)$.
\end{Proposition}

The equality here follows by counting 4-cliques minus an edge.

It immediately follows that the collinearity graph of a generalized quadrangle
(cf.~\cite{PayneThas84}) or partial quadrangle (cf.~\cite{Cameron74})
satisfies the 4-vertex condition (with $\alpha = \binom{\lambda}{2}$ and
$\beta = 0$). The same holds for a graph $\Gamma$ with $\lambda \le 1$.

If $\Gamma$ is locally strongly regular,
say with local parameters $(v',k',\lambda',\mu')$
(where clearly $v' = k$ and $k' = \lambda$),
then $\Gamma(x) \cap \Gamma(y)$ has valency $\lambda'$ (resp. $\mu'$)
when $x \adj y$ (resp. $x \nadj y$)
so that $\Gamma$ satisfies the 4-vertex condition with
$\alpha = \lambda\lambda'/2$ and $\beta = \mu\mu'/2$.

\subsection{A few rank 4 examples}
Below we give a small table with the parameters of some edge-transitive
rank 4 graphs satisfying the 4-vertex condition.
Except for the example with group $HJ.2$ due to Reichard \cite{Reichard00},
these do not seem to have been noticed in print.

{\medskip\footnotesize\noindent
\setlength{\tabcolsep}{5pt}
\begin{tabular}{@{\,}cccccccc@{~~~}l@{~~}l@{~~}l}
$v$ & $k$ & $\lambda$ & $\mu$ & $\lambda'$ & $\mu'$ & $\alpha$ & $\beta$ & group & name & ref \\
\hline
144 & 55 & 22 & 20 & - & 9 & 87 & 90 & ${\rm M}_{12} . 2$ \\
280 & 36 & 8 & 4 & - & 2 & 1 & 4 & ${\rm HJ} . 2$ &  & \cite{Reichard00} \\
300 & 104 & 28 & 40 & - & 8 & 78 & 160 & ${\rm PGO}_5(5)$ & $NO_5^-(5)$ & \S\ref{sec:Oex} \\
325 & 144 & 68 & 60 & - & 30 & 1153 & 900 & ${\rm PGO}_5(5)$ & $NO_5^+(5)$ & \S\ref{sec:Oex} \\
512 & 196 & 60 & 84 & 14 & 20 & 420 & 840 & $2^9 . {\rm \Gamma L}_3(8)$
    & dual hyperoval & \S\ref{sec:hyperovals} \\
729 & 112 & 1 & 20 & 0 & 0 & 0 & 0 & $3^6 . 2 . {\rm L}_3(4) . 2$ & Games graph & \cite{BrouwerVanLint84} \\
1120 & 729 & 468 & 486 & 297 & 306 & 69498 & 74358 & ${\rm PSp}_6(3) . 2$
     & disj.~t.i.~planes & \S\ref{sec:Sp6} \\
1849 & 462 & 131 & 110 & - & - & 2980 & 1845 & $43^2 {:} (42 {\times} {\rm D}_{22})$
     & power~diff.~set & \S\ref{sec:cyclo} \\
\end{tabular}\par}

\medskip
The numbers $\lambda',\mu'$ give the valency of the $\lambda$- and $\mu$-graphs
in case these are regular (and then $\alpha = \lambda\lambda'/2$ and
$\beta = \mu\mu'/2$).

\medskip
The examples on 144 and 729 vertices also satisfy the 5-vertex condition.

\subsection{Strongly regular graphs with strongly regular subconstituents}%
\label{srgsubs}
As we saw, graphs that are locally strongly regular satisfy
the 4-vertex condition.
Sometimes it follows that also the 2nd subconstituents must be
strongly regular.

\begin{Lemma} \label{lem:paras_for_local_srg}
Suppose that a strongly regular graph with parameters
$(v,k,\lambda,\mu) = (4t^2, 2t^2-\eps t,\maysplit t^2-\eps t, t^2-\eps t)$
(where $\eps = \pm1$)
has first subconstituents that are strongly regular with parameters
$(v',k',\lambda',\mu') =\big (2t^2-\eps t, t^2-\eps t,
\frac12 t(t - \eps), t(\frac12 t - \eps)\big)$.
Then its second subconstituents are strongly regular
with parameters $(v'',k'',\maysplit \lambda'',\mu'') = \big(2t^2 + \eps t - 1, t^2,
\frac12 t(t - \eps), \frac12 t^2\big)$.
\end{Lemma}

{\footnotesize
More generally, the spectrum of the 2nd subconstituent at any vertex
of a strongly regular graph follows from that of the 1st subconstituent%
---see \cite{CameronGoethalsSeidel78}, Theorem 5.1.\par\medskip}


Call the three parameter sets in the above lemma $A(\eps t)$, $B(\eps t)$,
and $C(\eps t)$, respectively. They occur again in \S\ref{binaryvs}.
The parameter sets $A(t)$ and $A(-t)$ are known as ({\em negative})
{\em Latin square parameters} ${\rm LS}_t(2t)$ (resp. ${\rm NL}_t(2t)$).
The complementary graphs have parameters ${\rm LS}_{t+1}(2t)$
(resp. ${\rm NL}_{t-1}(2t)$).

\medskip
Cameron, Goethals \& Seidel \cite{CameronGoethalsSeidel78}
studied the situation of a primitive strongly regular graph such that,
for some vertex, both subconstituents are strongly regular, and found
that such a graph either has a vanishing Krein parameter
$q_{11}^1$ or $q_{22}^2$, or has Latin square or negative Latin square
parameters. They conjectured that every non-grid example of the latter
has parameters as in the above lemma or has a complement with
these parameters.

\section{Survey of the known examples and results} \label{sec:knownknowns}
\subsection{Complements}
A graph satisfies the $t$-vertex condition if and only if its complement does.

\subsection{Generalized quadrangles}\label{sec:GQ}
Higman \cite{Higman71} observed that the collinearity graphs of
generalized quadrangles satisfy the 4-vertex condition
(and there are many examples that are not rank~3, cf.~\cite{Kantor86}).

{\smallskip\footnotesize
More generally the 4-vertex condition holds for partial quadrangles.
For example, the Hill graph with parameters
$(v,k,\lambda,\mu) = (4096,234,2,14)$
(derived from the cap constructed in \cite{Hill73})
has a rank 10 group and satisfies the 4-vertex condition with
$\alpha=1$, $\beta=0$.
\par\medskip}

Reichard \cite{Reichard15} showed that the collinearity graphs of
generalized quadrangles satisfy the 5-vertex condition, and that the
collinearity graphs of generalized quadrangles ${\rm GQ}(s,s^2)$
satisfy the 7-vertex condition.

{\smallskip\footnotesize
More generally the 5-vertex condition holds for partial quadrangles.
\par}

\subsection{Binary vector spaces with a quadratic form}\label{binaryvs}
The first non-rank-3 graph satisfying the 5-vertex condition was
constructed by A. V. Ivanov \cite{Ivanov89}:  a strongly regular graph
$\Gamma_0$ whose subconstituents
$\Gamma_1, \Gamma_2$   satisfy  the 4-vertex condition.
The parameters are as follows.

{\smallskip\footnotesize\noindent
\begin{tabular}{c|cccccccl}
 & $v$ & $k$ & $\lambda$ & $\mu$ & $\alpha$ & $\beta$ & $|G|$ & remarks \\
\hline
$\Gamma_0$ & 256 & 120 & 56 & 56 & 784 & 672
           & $2^{20} \cdot 3^2 \cdot 5 \cdot 7$
           & rank 4: $1+120+120+15$ \\
$\Gamma_1$ & 120 & 56 & 28 & 24 & 216 & 144
           & $2^{12} \cdot 3^2 \cdot 5 \cdot 7$
           & rank 4: $1+56+56+7$ \\
$\Gamma_2$ & 135 & 64 & 28 & 32 & 168 & 192
           & $2^{12} \cdot 3^2 \cdot 5 \cdot 7$
           & intransitive: $120+15$
\end{tabular}\par}

\medskip
In \cite{BrouwerIvanovKlin89} an infinite family of graphs $\Gamma^{(m)}$
($m \ge 1$) is constructed by taking as vertex set $\Ffss{2}{2m}$,
where vectors are adjacent when the line joining them~meets
the hyperplane at infinity in a fixed hyperbolic quadric
minus a maximal t.i.~subspace.
The graphs $\Gamma^{(m)}$ have parameters $A(2^{m-1})$ (see \S\ref{srgsubs}).
They have a rank 4 group (for $m \ge 4$) and satisfy the 4-vertex condition.

The local graphs $\Delta^{(m)}$ are strongly regular with parameters
$B(2^{m-1})$.
They have a rank 4 group (for $m \ge 4$) and satisfy the 4-vertex condition.

By Lemma \ref{lem:paras_for_local_srg} also the 2nd subconstituents
$\Eps^{(m)}$ are strongly regular, with parameters $C(2^{m-1})$.

We checked by computer that the graph $\Gamma^{(4)}$ is isomorphic to the
above $\Gamma_0$.

{\medskip\footnotesize
In \cite{Reichard00} it is shown that the graphs $\Gamma^{(m)}$
satisfy the 5-vertex condition.

In \cite{PechPech19} it is shown that the graphs $\Gamma^{(m)}$
are triplewise 5-regular, a.k.a. (3,5)-regular, where
$(s,t)$-regularity is the analog of the $t$-vertex condition
where $s$ instead of two vertices are distinguished.
It follows that the 2nd subconstituents $\Eps^{(m)}$
of the graphs $\Gamma^{(m)}$ also satisfy the 4-vertex condition.
\par}

\medskip
In \cite{Ivanov94}, two infinite families of graphs are constructed.
One is the above $\Gamma^{(m)}$. The second family has members
$\Sigma^{(m)}$ with vertex set $\Ffss{2}{2m}$, where vectors are
adjacent when the line joining them hits the hyperplane at infinity
either in a fixed elliptic quadric minus a maximal t.i.~subspace $S$
or in $S^\perp\backslash S$. The graphs $\Sigma^{(m)}$ have 
parameters $A(-2^{m-1})$,  have rank 5 (for $m \ge 5$),
%
%
and satisfy  the 4-vertex condition.
%
%

{\smallskip\footnotesize
Let $\Gamma(V,X)$ be the graph on a vector space $V$ where two vectors
are adjacent precisely when the joining line hits the subset $X$ of the
hyperplane $PV$ at infinity. Since $\Gamma(V,X)$ is strongly regular
if and only if $X$ is a 2-character set (\cite{Delsarte72a}), 
that is, if and only if $|X \cap H|$ takes only two distinct values
when $H$ runs through the hyperplanes of $PV$,
the set $(Q \setminus S) \cup (S^\perp \setminus S)$
is a 2-character set when $Q$ is an elliptic quadric,
and $S$ a maximal t.i. subspace.
%
%

Let $V$ be a  vector space over $\Ff_2$.
Then the local graph of $\Gamma(V,X)$ is the collinearity graph
of the partial linear space with point set $X$ and whose lines are the
projective lines (of size 3) contained in $X$.
%
\par\smallskip}

The local graphs $\Tau^{(m)}$ are strongly regular with parameters
$B(-2^{m-1})$. They are intransitive (for $m \geq 5$).
%
%

\mysqueeze{0.2pt}{It follows from Lemma \ref{lem:paras_for_local_srg} that
also the 2nd subconstituents $\Upsilon^{(m)}$
are~strongly} regular, with parameters $C(-2^{m-1})$.
There is a tower of graphs here:
If $\Upsilon$ is the 2nd subconstituent of $\Sigma^{(m)}$ at a vertex $x$,
and $s \in S$, then the local graph of $\Upsilon$ at its vertex $x+s$
is isomorphic to $\Sigma^{(m-1)}$. (For a proof, see Appendix A.)

{\smallskip\footnotesize
In \cite{Ivanov94} it is conjectured that the graphs $\Sigma^{(m)}$
satisfy the 5-vertex condition, and that the graphs $\Tau^{(m)}$
and $\Upsilon^{(m)}$ satisfy the 4-vertex condition.
The former was proved in \cite{Reichard00}.
The latter is proved in Appendix~A.
In \cite{PechPech19} it is announced that $\Sigma^{(m)}$ is even
$(3,5)$-regular, but we are not aware of a proof in print.
\par}

\subsection{Block graphs of Steiner triple systems}
Higman \cite{Higman71} investigated for which $v$-point Steiner triple systems
the block graph satisfies the 4-vertex condition.
He found that either the system is a projective space 
${\rm PG}(m,2)$ or $v$ is one of
9, 13, 25. In \cite{Kaski-et-al12} the cases 13 and 25 are ruled out,
so that the only other example is the affine plane ${\rm AG}(2,3)$.
The examples are rank 3.

\subsection{Smallest example}
In \cite{Klin-et-al05} it is shown that the smallest non-rank-3
strongly regular graphs satisfying the 4-vertex condition have
$v = 36$ vertices. There are three examples. All have
$(v,k,\lambda,\mu) = (36,14,4,6)$ and $\alpha=0$, $\beta=4$.

\subsection{Cyclotomic examples} \label{sec:cyclo}
Given $(q,e,J)$, where $e \,|\, (q-1)/2$ and $J$ is a set of
nonnegative integers, and a fixed primitive element $\eta$ of $\Ff_q$,
consider the cyclotomic graph with vertex set $\Ff_q$, where two elements
are adjacent when their difference is in
$D = \{ \eta^{ie+j} \mid 0 \le i < (q-1)/e, ~~ j \in J \}$.
In some cases this yields a strongly regular graph that satisfies
the 4-vertex condition. We give a few examples.
The examples on $11^2$ and $23^2$ vertices are due to
Klin \& Pech \cite{KlinPech}.

{\smallskip\footnotesize
\begin{tabular}{cccccccc}
$q$ & $p^f$ & $e$ & $J$ & $\eta$ & $\alpha$ & $\beta$ & rk \\
\hline
1849 & $43^2$ & 4 & $\{0\}$ & any & 2980 & 1845 & 4 \\
146689 & $383^2$ & 4 & $\{0\}$ & any & 11353825 & 10662960 & 4 \\
121 & $11^2$ & 6 & $\{0,1,2\}$ & any & 200 & 206 & 5 \\
625 & $5^4$ & 6 & $\{0,1,2\}$ & any & 5913 & 6022 & 5 \\
5041 & $71^2$ & 6 & $\{0,1,2\}$ & any & 395641 & 396270 & 5 \\
529 & $23^2$ & 8 & $\{0,1,2,3\}$ & $\eta^2 = \eta+4$ & 4215 & 4300 & 5 \\
\end{tabular}\par}

\smallskip
In all cases $q = p^f$ where $p$ is semiprimitive mod $e$
(that is, $e \,|\, (p^i+1)$ for some $i$),
so that the parameters of the strongly regular graph can be found
in \cite[Thm.~7.3.2]{BrouwerVM21}.

\section{Graphs from hyperovals}\label{sec:hyperovals}
In \cite{HuangHuangLin09}, Huang, Huang \& Lin constructed
various families of graphs. The complement of one of them 
can be described as follows (\cite{Brouwer16}). 
For $q = 2^m$, take $\Ffss{q}{3}$ as the vertex set of $\Gamma$.
Let $\pi$ be the plane at infinity of $\Ffss{q}{3}$. 
Let $H^*$ be a~dual hyperoval of $\pi$ (that is, a set of 
$q+2$ lines, no three on a point).
The plane $\pi$ is partitioned into two parts,
$\frac12 (q+1)(q+2)$ points on two lines of $H^*$
and $\frac12 q(q-1)$ exterior points on no line of $H^*$.
Two vertices of $\Gamma$ are adjacent when the  line joining them hits
$\pi$ in one of the   exterior points.
Then $\Gamma$ is strongly regular and has parameters
$$
(v, k, \lambda, \mu) = 
\big(q^3, \tfrac12 q(q-1)^2, \tfrac14 q(q-2)(q-3), \tfrac14 q(q-1)(q-2)\big).
$$
Its local graphs are strongly regular with parameters
$$
\big(\tfrac12 q(q-1)^2, \tfrac14 q(q-2)(q-3), \tfrac18 q(q^2-9q+22),
\tfrac18q(q-3)(q-4)\big).
$$
Hence, as noted in Section  \ref{sec:def}, $\Gamma$ satisfies the
$4$-vertex condition.
If $m=3$, then $\Gamma$ has rank $4$. 

\section{Disjoint t.i.~planes in symplectic 6-space}\label{sec:Sp6}
Let $V$ be a 6-dimensional vector space over $\Ff_q$,
provided with a nondegenerate symplectic form.
Let $\Gamma$ be the graph with as vertices the totally isotropic planes,
adjacent when disjoint.

\begin{Proposition}\label{disj-planes}
The graph $\Gamma$ is strongly regular, with parameters
$v = (q^3+1)(q^2+1)(q+1)$, $k = q^6$,
$\lambda = q^2(q^3-1)(q-1)$, $\mu = (q-1)q^5$.
If $q$ is even, then $\Gamma$ is rank $3$, otherwise rank $4$.
Its local graph $\Delta$ is strongly regular with parameters
$v' = k$, $k' = \lambda$, $\lambda' = \mu' - q^2(q-2)$
and $\mu' = q^2(q-1)(q^3-q^2-1)$.
It follows that $\Gamma$ satisfies the $4$-vertex condition.
\end{Proposition}

For convenience, we give the parameters of $\bar{\Delta}$,
the complement of $\Delta$:\\
$\bar{v} = q^6$, $\bar{k} = (q^2+1)(q^3-1)$,
$\bar{\lambda} = q^4+q^3-q^2-2$, $\bar{\mu} = q^4+q^2$.

{\footnotesize\medskip
\Proof
The dual polar graph $\Sigma$ belonging to ${\rm Sp}_6(q)$
is distance-regular of diameter 3 and has eigenvalue $-1$.
%
%
It follows that its distance-3 graph $\Gamma$ is strongly regular
(see \cite{BCN}, Prop.~4.2.17).
More generally, the distance 1-or-2 graph of the symplectic 
dual polar space ${\rm Sp}_{2m}(q)$ is distance-regular
(cf.~\cite{BCN}, Prop.~9.4.10).
For $m=3$ it is the complement of $\Gamma$.
%
%
%
%

For any vertex $x$, the subgraph induced by $\Sigma$ on $\Sigma_3(x)$
is isomorphic to the symmetric bilinear forms graph on $\Ffss{q}{3}$
(see \cite{BCN}, Prop.~9.5.10).
If $q$ is odd, then distance $j$ ($j=0,1,2,3$) in $\Sigma_3(x)$
corresponds to ${\rm rk}(f-g) = j$  in the symmetric bilinear forms graph
and hence to distance $\lfloor (j+1)/2 \rfloor$ in the quadratic forms
graph (see \cite{BCN}, \S9.6). It follows that $\Delta$ is the
complement of the quadratic forms graph, and has parameters as claimed.


If $q$ is even, then $\Gamma$ is rank 3 (by triality, it is
the complement of the $O_8^+(q)$ polar graph), and $\Delta$
is the complement of the rank 3 graph $\smash{VO_6^+(q)}$,
with parameters as claimed.
%
\qed
\par}

\section{Nonsingular points joined by a tangent}\label{sec:Oex}
Let $V$ be a vector space of dimension $2m+1$ over  $\Ff_q$ 
with $q$   odd, and let $Q$ be a nondegenerate quadratic form on $V$.
We  also use $Q$ as the symbol for the set of singular projective points.

The projective space $PV$ has $(q^{2m+1}-1)/(q-1)$ points,
$(q^{2m}-1)/(q-1)$ singular, and $q^{2m}$ nonsingular.
The nonsingular points come in two types: 
there are $\frac12 q^m(q^m + \eps)$ points of type $\eps$
(where $\eps = \pm 1$), with $\eps=+1$ (resp. $-1$) for points $x$
for which $x^\perp$, the hyperplane of points orthogonal to $x$,
is hyperbolic (resp. elliptic).

Consider the graph $NO_{2m+1}^\eps(q)$ that has as vertex set the set
of nonsingular points of type $\eps$, where two points are adjacent
when the joining line is a tangent.

\begin{Proposition} {\rm (Wilbrink \cite{Wilbrink},
cf.~\cite{BrouwerVanLint84})}
Let $m \ge 2$. The graph $NO_{2m+1}^\eps(q)$ is strongly regular
with parameters
$v = \frac12 q^m(q^m + \eps)$,
$k = (q^{m-1}+\eps)(q^m-\eps)$,
$\lambda = 2(q^{2m-2}-1) + \eps q^{m-1}(q-1)$,
$\mu = 2q^{m-1}(q^{m-1}+\eps)$.
\end{Proposition}

\vspace{-0.2cm}
For $m=1$, $\eps=-1$ the graph is edgeless.
For $m=1$, $\eps=1$ we have the triangular graph $T(q+1)$.
Wilbrink also handled the case of even $q$.
We give an explicit proof here;
for a different and more general proof see \cite{BannaiHaoSong90}.

{\medskip\footnotesize
\Proof
The neighbors of a vertex $x$ lie on the tangents
joining $x$ with a singular point of $x^\perp$,
and $x^\perp$ has $(q^{m-1}+\eps)(q^m-\eps)/(q-1)$ singular points.
This gives the value of $k$.

A common neighbor $z$ of two adjacent vertices $x,y$
lies on the line $xy$ (and there are $q-2$ choices) or on
some other tangent $T$ on $x$.
In the latter case the plane $\< x,y,z \>$ meets $Q$
in a conic or double line. If it is a conic, then $z$ is
uniquely determined on $T$ by the fact that $yz$ is the tangent
on $y$ other than $xy$. If it is a double line, then each
nonsingular point of $T \setminus \{x\}$ is suitable.
Let $p$ be the singular point on $xy$.
Then $\{p,x\}^\perp/\< p \>$ is a nondegenerate
$(2m-2)$-space of type $\eps$, and has
$a = (q^{m-2}+\eps)(q^{m-1}-\eps)/(q-1)$ singular points.
It follows that $xy$ is in $a$ planes that hit $Q$ in a double line,
and in $q^{2m-2}$ planes that hit $Q$ in a conic.
Consequently, $\lambda = q-2 + q^{2m-2} + (q-1)qa$,
as desired.

A common neighbor $z$ of two nonadjacent vertices $x,y$
determines a nondegenerate plane $\pi = \< x,y,z \>$ in which
$xz$ and $yz$ are tangents, so that $x,y,z$ are exterior points.
Now $x,y$ are on two tangents each, and $\pi$ contains 4
common neighbors of $x,y$.
If $Q$ is a quadratic form on a $(2m+1)$-space, then a point $p$
is exterior if and only if $(-1)^m \det(Q)\,Q(p)$ is a nonzero square.
In order to have $p$ exterior in $\pi$ but a $\eps$-point in $V$,
the $(2m-2)$-space $\pi^\perp$ must be an $\eps$-subspace of the
$(2m-1)$-space $\{x,y\}^\perp$. Since there are
$b = \frac12 q^{m-1}(q^{m-1}+\eps)$ such $\eps$-subspaces,
we find $\mu = 4b$, as desired. \qed
\par}

\medskip
The automorphism group ${\rm P\Gamma{}O}_{2m+1}(q)$ 
of the graph contains ${\rm PGO}_{2m+1}(q)$.
The latter has $(q+3)/2$ orbits on pairs of vertices \cite{BannaiHaoSong90}.
Hence, the graph has rank $(q+3)/2$ if $q$ is prime.

For $m=2, \eps = -1$, this is the collinearity graph of a semi-partial geometry
found by Metz. Its lines have size $s+1 = q$ and there are $t+1 = q^2+1$
lines on each point. Each point outside a line has either $0$ or
$\alpha = 2$ neighbors on the line.
See Debroey \cite{Debroey78}, voorbeeld 1.1.3d,
and Debroey-Thas \cite{DebroeyThas78}, example 1.4d,
and Hirschfeld-Thas \cite{HirschfeldThas80}, p.~268,
and Brouwer-van\,Lint \cite{BrouwerVanLint84}, \S7A,
and Brouwer-Van\,Maldeghem \S8.7, example (ix).

For $m=2, \eps = +1$ this is the collinearity graph of a geometry
with $t+1 = (q+1)^2$ lines of size $s+1 = q$ on each point,
such that each point outside a line has 0, 2, or $q$ neighbors on the line 
(\cite{BrouwerVanLint84}, \S7B).

\medskip
We shall prove that these graphs satisfy the 4-vertex condition.
First a lemma.

\begin{Lemma}\label{lem:tri_adj_sd}
    Let $S$ be a solid such that $Q\restrictedto{S}$ is nondegenerate.
    Let $x, y, z$ be distinct nonsingular points of the same type $\eps$
    such that $\< z, x \>$ and $\< z, y \>$ are tangents
    and $\< x, y \>$ is nondegenerate.
    Put $\pi = \< x,y,z \>$. Then there are 
    either $0$ or $2$ nonsingular points $w \in S \setminus \pi$ of type $\eps$
    such that $\< x, w \>$, $\< y, w \>$, 
    and $\< z, w \>$ are tangents.
    For $x, y, z$ given, the number of $w$ only depends on the type of $S$.
    It equals $2$ if and only if the nonzero number
    $2(\frac{B(z,z)B(x,y)}{B(x,z)B(y,z)}-1) \det (Q\restrictedto{S})$
    is a square.
\end{Lemma}

\Proof 

Replace $x$ by $\frac{B(z,z)}{B(x,z)} x$ and
$y$ by $\frac{B(z,z)}{B(y,z)} y$.
Then $B(x,z) = B(z,z) = B(y,z)$.
Put $x_0 = x-z$, $y_0 = y-z$, $w_0 = w-z$, then
$B(x_0,z) = B(y_0,z) = B(w_0,z) = 0$.
Since the lines $\< z, x \>$, $\< z, y \>$, and
$\< z, w \>$ are tangents, the points $x_0,y_0,z_0$ are singular,
that is, $Q(x_0) = Q(y_0) = Q(w_0) = 0$.
The line $\< x, w \>$ is a tangent, so $Q(x+tw) = 0$
has a unique solution $t$. Now
\begin{align*}
  Q(x+tw) &= Q(z + x_0 + t(z + w_0)) = Q((1+t)z + x_0 + t w_0)\\
  &= (1+t)^2 Q(z) + Q(x_0 + t w_0) = (1+t)^2 Q(z) + tB(x_0, w_0).
\end{align*}
It follows that $(2+\frac{B(x_0, w_0)}{Q(z)})^2 = 4$,
that is $\frac{B(x_0, w_0)}{Q(z)} \in \{ 0, -4 \}$.

As $Q\restrictedto{S}$ is nondegenerate, $z^\perp \cap S$ is
a nondegenerate plane.
If $B(x_0, w_0) = 0$, then $\< x_0, w_0 \>$ is a
totally singular line in this plane, impossible.
Hence, $B(x_0, w_0) = -4 Q(z)$. Similarly, $B(y_0, w_0) = -4 Q(z)$.

In the plane $z^\perp \cap S$, let $u$ be the point of intersection of
the tangents through the points $x_0$ and $y_0$ and write
$w_0 = ax_0 + by_0 + cu$. Then $B(x_0, u) = B(y_0, u) = 0$ and
$-4 Q(z) = B(x_0, w_0) = B(x_0, ax_0 + by_0 + cu) = bB(x_0, y_0)$.
Similarly, $-4 Q(z) = B(y_0, w_0) = aB(x_0, y_0)$,
so that $a = b = \frac{-4Q(z)}{B(x_0, y_0)}$, independent of $w$. Also, 
\begin{align*}
    0 &= Q(w_0) = Q(ax_0 + by_0 + cu) = abB(x_0, y_0) + c^2 Q(u)
= \frac{16Q(z)^2}{B(x_0, y_0)} + c^2 Q(u).
\end{align*}
If $-B(x_0, y_0)Q(u)$ is a square, then we have two solutions for $c$
(so also $w_0$ and, therefore, $w$) and otherwise none.
Since $u$ is an exterior point in the plane $\sigma = z^\perp \cap S$,
the number $-Q(u) \det Q\restrictedto{\sigma}$ is a square.
Also, $\det Q\restrictedto{S} = Q(z) \det Q\restrictedto{\sigma}$
and $B(x,y) = B(x_0,y_0) + B(z,z)$.
\qed

\begin{Proposition}
The graph $NO_{2m+1}^\eps(q)$ satisfies the 4-vertex condition.
\end{Proposition}

\Proof
By Proposition \ref{prop:sims} it suffices to check for $x \ne y$ that
the number of edges in $\Gamma(x) \cap \Gamma(y)$ does not depend on
the choice of the points $x,y$, but only on whether $x,y$
are adjacent or not.

Since ${\rm Aut}~\Gamma$ is edge-transitive,
we only need to check $\Gamma(x) \cap \Gamma(y)$ for $x \nadj y$.

\smallskip 
Claim: this subgraph $\Gamma(x) \cap \Gamma(y)$
is regular of valency $4q^{2m-3} + 3\eps q^{m-1} - 4\eps q^{m-2} - 1$.
In other words, this is the value of $\mu$ in the local graph
(which is regular, but not strongly regular).

If $x \adj z \adj y$, $x \nadj y$, then $\pi = \< x,y,z \>$
is a nondegenerate plane in which the common neighbors of $x,y$
form a 4-cycle, so that $x,y,z$ have two common neighbors in $\pi$,
say $a$ and $b$.

The plane $\pi$ lies in $(q^{2m-3}-\eps q^{m-2})/2$ solids of type
$O^-(4, q)$, equally many solids of type $O^+(4, q)$,
and $(q^{m-2}+\eps)(q^{m-1}-\eps)/(q-1)$ degenerate solids.

If $S$ is a degenerate solid through $\pi$ with apex $p$, we see that
$w \in S \setminus \pi$ is in $\Gamma(x) \cap \Gamma(y) \cap \Gamma(z)$
if and only if gets projected from $p$ onto an element of $\{ a, b, z \}$
in $\pi$. Hence,
$|\Gamma(x) \cap \Gamma(y) \cap \Gamma(z) \cap S \setminus \pi| = 3(q-1)$.
Hence, the total number of choices for $w$ equals $3(q^{m-2}+\eps)(q^{m-1}-\eps)$.

Now let $S$ be a nondegenerate solid on $\pi$,
and let $p = S \cap \pi^\perp$.
By Lemma \ref{lem:tri_adj_sd}, the number of $w$ in $S$ is 0 or 2,
depending on the determinant of $Q$ restricted to $S$.
Since $\pi^\perp$ contains equally many points $p$ with $Q(p)$
a square as with $Q(p)$ a non-square, the total number of choices
for $w$ equals the number of choices for $p$ which is $q^{2m-3}-\eps q^{m-2}$.

So the induced subgraph on $\Gamma(x) \cap \Gamma(y)$ has valency
$2 + 3(q^{m-2}+\eps)(q^{m-1}-\eps) + (q^{2m-3}-\eps q^{m-2}) 
= 4q^{2m-3} + 3\eps q^{m-1} - 4\eps q^{m-2} - 1$.
\qed

%

\section{Polar switching}\label{sec:switched_col}
{\smallskip\footnotesize
A {\em polar space} is a partial linear space such that for each line $L$
any point outside $L$ is collinear to either all or precisely one of the
points of $L$. A {\em singular subspace} is a line-closed set of points,
any two of which are collinear. The polar space is called
{\em nondegenerate} when no point is collinear to all points.
Finite nondegenerate polar spaces are the sets of totally isotropic
(t.i.) or totally singular (t.s.) points and lines in a vector space
over a finite field provided with a suitable symplectic, quadratic or
hermitian form. The {\em rank} of the polar space is the
(vector space) dimension of its maximal singular subspaces.
\par}

\medskip
Let $\bP$ be a nondegenerate polar space of rank $d \ge 3$
in a vector space $V$ over $\Ff_q$. Its collinearity graph
$\Gamma_0$ is strongly regular and satisfies the 4-vertex condition
(since it is rank 3).
We shall construct cospectral graphs that satisfy the 4-vertex condition
(but are not rank 3) by a switching construction.
Let $x^\perp$ be the set of points collinear with $x$ (including $x$ itself).

{\medskip\footnotesize
Suppose $U$ is a maximal singular subspace of $\bP$
(i.e., a maximal clique in $\Gamma_0$), and let $H_1,H_2$ be two
hyperplanes of $U$. We can redefine adjacency and make the points
$x$ with $x^\perp \cap U = H_1$ or $H_2$ adjacent to the points in
$H_2$ or $H_1$, respectively, and leave all other adjacencies unchanged.
This is an example of WQH-switching (Wang, Qiu \& Hu \cite{WangQiuHu19},
cf.~\cite{Ihringer19}) and yields a graph cospectral with $\Gamma_0$.
One can repeat this interchange of hyperplanes and get arbitrary
permutations of all hyperplanes. We generalize this, even allowing
different designs on $U$.
\par}

\subsection{Construction}

Let $P$ be the point set of $\bP$, and let the subset
$U$ be (the set of points of) a totally isotropic $d$-space.
Let $\bD$ be a symmetric design with the same parameters 
as the symmetric design of points and hyperplanes of ${\rm PG}(d-1, q)$, so its
parameters are
2-$\big(\downstrut\frac{q^{d}-1}{q-1}, \frac{q^{d-1}-1}{q-1},
\frac{q^{d-2}-1}{q-1} \big)$. 
Let $\phi$ be a bijection from the set $\cH$ of hyperplanes of $U$
to the blocks of $\bD$. We assume that the points of $U$ are also
the points of $\bD$.

Following ideas in \cite{Kantor94} 
and \cite{DempwolffKantor08}
we define a graph $\Gamma_\phi$ on the vertex set of $\Gamma_0$ as follows:
\begin{enumerate}
 \item Vertices in $U$ are pairwise adjacent.
 \item Distinct vertices $x,y \notin U$ are adjacent if $x \in y^\perp$.
 \item Vertices $x \in U$, $y \notin U$ are adjacent if $x \in (y^\perp \cap U)^\phi$.
\end{enumerate}

Clearly, $\Gamma_\phi = \Gamma_0$ if we take the hyperplanes of $U$
for the blocks of $\bD$ and $\phi$ as the identity.

\begin{Theorem}
 The graph $\Gamma_\phi$ is strongly regular with
 the same parameters as the classical graph $\Gamma_0$.
\end{Theorem}
{\footnotesize
\Proof 
Let $x$ and $y$ be any two vertices. We show that
the number of common neighbors $z$ of $x,y$ in $\Gamma_\phi$
does not depend on $\phi$ (but depends on whether $x,y$ are equal,
adjacent or nonadjacent in $\Gamma_\phi$).

If $x,y \in U$, then any $z \in U$ is a common neighbor.
The number of $z \in P \setminus U$ such 
that $x,y \in (z^\perp \cap U)^\phi$ does not depend on $\phi$:
each hyperplane $H$ of $U$ such that $x, y \in H^\phi$
contributes $| H^\perp \setminus U |$ such $z$.

Suppose that $x,y \notin U$. Then we are counting the $z$ in
$(x^\perp \cap U)^\phi \cap (y^\perp \cap U)^\phi$,
and also the $z$ in $(x^\perp \cap y^\perp) \setminus U$.
The numbers of such $z$ does not depend on $\phi$.

The remainder of the proof concerns the case $x \in U$, $y \notin U$.
If $z \in U$  then the requirements are $z \ne x$ and
$z \in (y^\perp \cap U)^\phi$. The number of such $z$
does not depend on $\phi$.

So we need to count the  $z \notin U$.  First set $I := y^\perp \cap U$,
so $Y := \< y, I \>$ is totally isotropic.
If $z \in Y$  then $I^\phi = (z^\perp \cap U)^\phi$,
and $x,z$ are adjacent if and only if $x,y$ are adjacent. The number of
such $z$ is independent of $\phi$.

It remains to count the $z$ in $y^\perp \setminus Y$ such that
$x \in (z^\perp \cap U)^\phi$; here $z^\perp \cap U \ne I$ as $z\notin Y$.
Let $H \neq I$ be a hyperplane of $U$ such that $x \in H^\phi$.
The number of $H$ does not depend on $\phi$
(note that $x \in I^\phi$ if and only if $x,y$ are adjacent in $\Gamma_\phi$).
We show that the number of $z$ in $y^\perp \setminus Y$
with $z^\perp \cap U = H$ does not depend on $\phi$ or $H$.
Using bars to project $(H \cap I)^\perp$ into the nondegenerate
rank 2 polar space $(H \cap I)^\perp/(H \cap I)$, we see
totally isotropic lines $\bar U$ and $\bar Y$ meeting at a point $\bar I$,
and a nondegenerate 2-space $\< \bar y , \bar H \>$;
the number of $\bar z$  in  $\< \bar y , \bar H \>^\perp\backslash \bar I$
does not depend on $\phi$ or $H$, so neither does the number of required $z$. 
\qed
\par\smallskip}

\subsection{Isomorphisms}

\subsubsection*{Emptying bijections $\phi$}
Call a vertex $e \in U$ {\em emptying} for $\phi$ if
$\bigcap \{ H \mid H \in \cH,~ e \in H^\phi \} = \emptyset$.
Call $\phi$ {\em emptying} if the subspace $U$ is spanned by 
emptying vertices.

Call a vertex $f \in U$ {\em dually emptying} for $\phi$ if
$\bigcap \{ H^\phi \mid f \in H \in \cH \} = \emptyset$.
Call $\phi$ {\em dually emptying} if the subspace $U$ is spanned by
dually emptying vertices.

\medskip
{\footnotesize
If $a$ is not emptying, then
$\bigcap \{ H \mid H \in \cH,~ a \in H^\phi \} = \{b\}$ for some vertex $b$.
If $b$ is not dually emptying, then
$\bigcap \{ H^\phi \mid b \in H \in \cH \} = \{a\}$ for some vertex $a$.
This establishes a 1-1 correspondence between not emptying vertices $a$
and not dually emptying vertices $b$.
\par}

\begin{Proposition}\label{prop:all_is_empty}
If a permutation $\phi$ of $\cH$ is not dually emptying,
then it is in ${\rm P\Gamma{}L}(U)$.
\end{Proposition}
{\footnotesize
\Proof
Let $E$ denote the set of emptying vertices of $U$, and put
$A = U \setminus E$. Let $F$ denote the set of dually emptying
vertices of $U$, and put $B = U \setminus F$.
Let $\psi \colon B \to A$ be the 1-1 correspondence found above.
%
We show that if $L$ is a line in $U$ with $|L \cap B| \ge q$,
then $L \subseteq B$ and $L^\psi$ is a line.

Indeed, let $b,b' \in L \cap B$ and set $M = \<b^\psi, b'^\psi\>$.
Then $L \subseteq H$ is equivalent to $M \subseteq H^\phi$
so that $(L \cap B)^\psi = M \cap A$.
If all points of $L$ are in $B$ with a single exception $w$,
then all points of $M$ are in $A$ with a single exception $v$,
and all hyperplanes $H$ with $w \in H$ satisfy $v \in H^\phi$
(since every line meets every hyperplane), and $v = w^\psi$,
that is, $w$ was no exception.

If $\phi$ is not dually emptying, then there
exists a hyperplane $H$ such that $U \setminus H \subseteq B$.
By the above this implies $B = U$ and $\psi$ is in
${\rm P\Gamma{}L}(U)$ and induces $\phi$ on the set $\cH$.
\qed
\par}

\subsubsection*{Large cliques}

We use the presence of maximal cliques of various sizes to study
the structure of the graphs $\Gamma_\phi$ when $\phi$ is a permutation.

\smallskip
Abbreviate the size $\frac{q^i-1}{q-1}$ of an $i$-space with $m_i$,
so that maximal singular subspaces have size $m_d$.
Since $m_d$ is the Delsarte-Hoffman upper bound for the size of cliques
in $\Gamma_\phi$, each vertex outside a clique of this size is
adjacent to precisely $m_{d-1}$ vertices inside,
cf. \cite[Proposition 1.1.7]{BrouwerVM21}.

\begin{Lemma}\label{maxcliques}
Let $d \ge 3$.

(i) If $M \ne U$ is a maximal singular subspace of $\bP$,
then $C = (M \setminus U) \,\cup\,
\bigcap \{ H^\phi \mid M \cap U \subseteq H \in \cH \}$
is a maximal clique in $\Gamma_\phi$ of size at least $q^{d-2}(q+1)$
(and $C \setminus U = M \setminus U$).

(ii) If $C \ne U$ is a maximal clique in $\Gamma_\phi$ of size
at least $q^{d-2}(q+1)$, then $M = \< C \setminus U \>$ is
a maximal singular subspace of $\bP$.

If, moreover, $|C| = m_d$, then $M \setminus U = C \setminus U$.

\end{Lemma}
\Proof
(i)
Let $M$ be a maximal singular subspace other than $U$.
Then $C = (M \setminus U) \cup
\bigcap \{ H^\phi \mid M \cap U \subseteq H \in \cH \}$
is the largest clique in $\Gamma_\phi$ containing $M \setminus U$.
(Indeed, the set of hyperplanes of $U$ of the form $m^\perp \cap U$
where $m \in M \setminus U$ equals the set of hyperplanes containing
$M \cap U$, so $C$ is a clique. No further point outside $U \cup C$ can be
adjacent to all of $C$, since $|M \setminus U| > m_{d-1}$.)
If $\dim M \cap U = d-1$, then $|C| = |M| = m_d$.
If $\dim M \cap U \le d-2$, then $|C| \ge |M \setminus U|
\ge m_d - m_{d-2} = q^{d-2}(q+1)$.

\smallskip
(ii)
Let $C \ne U$ be a maximal clique of size at least $q^{d-2}(q+1)$.
If $|C \setminus U| \le m_{d-1}$, then
$|C \cap U| \ge q^{d-2}(q+1) - m_{d-1} > m_{d-2}$.
The set $C \cap U$ is the intersection of sets $H^\phi$,
each of size $m_{d-1}$, and any two distinct such sets meet in
$m_{d-2}$ points. It follows that no two different $H$ occur, that is,
$H = c^\perp \cap U$ is independent of the choice of $c \in C \setminus U$.
Now $C$ is contained in, and hence equals, $H^\phi \cup (C \setminus U)$,
and $|C \setminus U| = m_d - m_{d-1} > m_{d-1}$, a contradiction.

If $S$ is a clique in $\Gamma_0$, then also $\< S \>$ is a clique
in $\Gamma_0$. In particular, $\< C \setminus U \>$ is
a singular subspace. It is maximal since $|\< C \setminus U \>| > m_{d-1}$.

\smallskip
If $|C| = m_d$, then each vertex outside $C$ is adjacent to
precisely $m_{d-1}$ vertices inside. Hence no point
outside $C \cup U$ can be adjacent to all of $C \setminus U$.
\qed

\begin{Lemma}\label{lem:detU}
If the permutation $\phi$ is dually emptying,
then $U$ is uniquely determined within the graph $\Gamma_\phi$.
\end{Lemma}
\Proof
The subspace $U$ is a clique of size $m_d$ in $\Gamma_\phi$,
with the two properties

(i) in the subgraph induced on its complement $P \setminus U$
all maximal cliques $N$ have size $m_d - m_i$
(where $m_i = | \< N \> \cap U |$) for some $i$, $0 \le i \le d-1$,
and

(ii) the number of maximal cliques of size $m_d$ disjoint from $U$
equals the number of maximal singular subspaces disjoint from any
given one.

\smallskip
Let $E \ne U$ be a clique of $\Gamma_\phi$ of size $m_d$ 
with the same two properties. 
First we use (i) to see that $E \cap U$ must be a hyperplane in $U$.

\smallskip
Since $E$ is a maximal clique, and $\phi$ is
a permutation, $E \cap U$ is an intersection of hyperplanes
and hence a subspace of $U$. By hypothesis, we can find a 
dually emptying point $f$ of $U$ not in $E$.
If $g \in f^\perp \cap (E \setminus U)$
($g$ will exist unless $f^\perp \cap E = U \cap E$)
and $M$ is a maximal singular subspace containing $f$ and $g$,
and meeting $U$ in $\{f\}$, then $C = M \setminus \{ f \}$
is a maximal clique in $\Gamma_\phi$ of size $m_d-1$.
And $N = C \setminus E$ is a maximal clique in $P \setminus E$
of size $m_d - m_i - 1$ in case $|M \cap E| = m_i$.
(Note that $C \setminus U = M \setminus U$.)

\smallskip
{\footnotesize
Why is N maximal? No point can be added since $|N| > m_{d-1}$,
unless $q=2$ and $|N|=|M \cap E|=m_{d-1}$. In that case,
no point outside U can be added since $\<N\>=M$.
And no point inside $U$ can be added since $N$ determines all
hyperplanes on $f$, and $f$ is dually emptying.
\par}

\smallskip
Since $M \cap E \ne \emptyset$, we have $1 \le i \le d-1$,
and $m_d - m_i - 1$ is not of the form $m_d - m_h$, violating (i).
Therefore, $f^\perp \cap E = U \cap E$, so that
$H = \<E \setminus U \> \cap U$ and $H^\phi = E \cap U$ are hyperplanes.

\smallskip
Now we use (ii) to arrive at a contradiction.

We claim that if a maximal clique $F$ of size $m_d$ is disjoint from
$E$, then $\< F \setminus U \>$ is disjoint from $\< E \setminus U \>$.
Suppose not. Since $\< E \setminus U \> \setminus U = E \setminus U$
and $\< F \setminus U \> \setminus U = F \setminus U$ by Lemma \ref{maxcliques}(ii),
a common vertex must lie in $U$.
If $\< F \setminus U \>$ meets $U$ in $m_e$ vertices with $e \geq 2$, 
then $F$ meets $U$ in a subspace of dimension $e$, but that would meet
$H^\phi$, impossible.
So, $\< F \setminus U \>$ meets $U$ in a singleton $\{ f \}$
on the hyperplane $H$.
As $F$ has size $m_d$, $f$ is not dually emptying,
so $\bigcap \{ H^\phi \mid f \in H \} = \{ f' \}$ for some point $f'$.
Now $f' \in E \cap F$, a contradiction. This shows our claim.

By the claim and Lemma \ref{maxcliques}, we have an injection from the set
of maximal cliques of size $m_d$ disjoint from $E$ into the set of
maximal singular subspaces disjoint from $\< E \setminus U \>$.
Since $E$ satisfies (ii), both sets have the same size, so
the injection is also a surjection.

On the other hand, since $\phi$ is dually emptying, there is a
dually emptying point $o$ in $U \setminus H$. This $o$ lies in
a maximal singular subspace $O$ disjoint from $\< E \setminus U \>$,
and this $O$ is not in the image of the surjection. Contradiction.
\qed

\begin{Lemma}\label{lem:extending}
Let $\bP$ be a nondegenerate polar space with point set $P$,
and $U$ a maximal totally isotropic subspace.
Let $h \colon P \setminus U \to P \setminus U$ be a bijection preserving
collinearity. Then $h$ can be uniquely extended to an automorphism
$h'$ of $\bP$.
\end{Lemma}
{\footnotesize
\Proof
Indeed, we can extend $h$ as follows.
For $u \in U$, let $R$ be a maximal t.i.~subspace with $U \cap R = \{u\}$.
Then $R \setminus \{u\}$ is a subspace of $\bL$ of size
$|U| - 1$ and is mapped by $h$ to a similar subspace $S$.
In $\bP$ this subspace is contained in a unique maximal t.i.~subspace
$\<S\>$ ($= S^\perp$) and we can define $h'(u) = v$ when
$\<S\> \setminus S = \{v\}$.

This is well-defined: if $R'$ is a maximal t.i.~subspace
with $U \cap R' = \{u\}$ and $R$, $R'$ meet in codimension 1,
and $h$ maps $R' \setminus \{u\}$ to $S'$, then
$\< S \cap S' \> = (S \cap S') \cup \{v\}$.
Since the graph on such subspaces $R$, adjacent when they meet
in codimension 1, is connected, $v$ is well-defined.

This preserves orthogonality: if $u \in x^\perp$, then there is
a maximal t.i.~subspace $R$ containing $u,x$ with $R \cap U = \{u\}$.
Now $h(u) = v$ lies in the t.i.~subspace $\< h(R \setminus \{u\}) \>$
which also contains $h(x)$.
\qed
\par}

\begin{Proposition}\label{prop:isos_kantor}
  \mysqueeze{0.378pt}{Let $\bP$ be a nondegenerate polar space and $U$ a maximal t.i.~subspace.}
  Let $\phi$ and $\chi$ be permutations of $\cH$ such that
  $\Gamma_\phi$ is isomorphic to $\Gamma_\chi$. Then $\phi$ and $\chi$ are 
  in the same ${\rm P\Gamma{}L}(U)$-double coset in Sym$(\cH)$.
\end{Proposition}


\Proof
If $\phi \in {\rm P\Gamma{}L}(U)$, then $\Gamma_\phi$ is isomorphic to
$\Gamma_0$ and its group of automorphisms is transitive on the set of
maximal singular subspaces.
If $\phi \notin {\rm P\Gamma{}L}(U)$, then according to
Lemma \ref{lem:detU} and Proposition \ref{prop:all_is_empty}
the maximal singular subspace $U$ can be recognized
in $\Gamma_\phi$, and hence $\Gamma_\phi$ is not isomorphic to $\Gamma_0$.
Since by assumption $\Gamma_\phi$ and $\Gamma_\chi$ are isomorphic,
either both or neither are isomorphic to $\Gamma_0$.
In the former case both $\phi$ and $\chi$ are in ${\rm P\Gamma{}L}(U)$
and the claim holds. Assume in the following that $\phi$ and $\chi$
are not in ${\rm P\Gamma{}L}(U)$.

\medskip
We have the set $P$, the point set of $\bP$,
with three structures defined on it.
The polar space structure $\bP$, with relation $\perp$,
and the two graph structures $\Gamma_\phi$ and $\Gamma_\chi$.
We translate what it means for $\Gamma_\phi$ and $\Gamma_\chi$
to be isomorphic in terms of the polar space.

Let $g: \Gamma_\phi \rightarrow \Gamma_\chi$ be an isomorphism.
By Lemma \ref{lem:detU}, it sends $U$ to itself.

The number of common neighbors of a triple of points in $U$
equals $\lambda-1$ for collinear triples and is smaller
for noncollinear triples. It follows that $g$ preserves
projective lines in $U$, and hence induces a permutation
$\bar{g}$ of $\cH$ that is in ${\rm P\Gamma{}L}(U)$.

Let $h$ denote the restriction of $g$ to $P \setminus U$.
Then $h$ preserves collinearity (since we have
$\{x,y,z\}^\perp \cap (P \setminus U) = \{x,y\}^\perp \cap (P \setminus U)$
for a triple of pairwise adjacent points $x,y,z$ of $P \setminus U$
if and only if $x,y,z$ are collinear).
By Lemma \ref{lem:extending}, $h$ can be uniquely extended
to an automorphism $h'$ of $\bP$.

\smallskip
Let $\bar{h}$ be the permutation of $\cH$ induced by $h'$.
Then $\bar{h} \in {\rm P\Gamma{}L}(U)$.

For $x \in U$ and $y \notin U$, if $x$ and $y$ are adjacent in
$\Gamma_\phi$, then $x^g$ and $y^g$ are adjacent in $\Gamma_\chi$.
This says that $x \in (y^\perp \cap U)^\phi$ implies that
$x^g \in (y^{g\perp} \cap U)^\chi$: $g$ maps the points of any hyperplane
of $U$ to the points of another hyperplane.
Then $(y^\perp \cap U)^{\phi g} = (y^{g\perp} \cap U)^\chi
 = (y^{h\perp} \cap U)^\chi = (y^\perp \cap U)^{\bar{h}\chi}$,
 so that $\phi\bar{g} = \bar{h}\chi$.
\qed

\begin{Theorem}\label{thm:low_bnd}
Let $d \ge 3$. There are at least $q^{d-2}!$ pairwise 
nonisomorphic strongly regular graphs having the same parameters
as the collinearity graph $\Gamma_0$ of the polar space $\bP$.
\end{Theorem}

\Proof
Let $q = p^e$, where $p$ is prime.
Then $|{\rm P\Gamma{}L}(U)| < e q^{d^2}$.
In view of Proposition \ref{prop:isos_kantor}, we have obtained at least
$m_d!/|{\rm P\Gamma{}L}(U)|^2 > q^{d-2}!$
pairwise nonisomorphic strongly regular graphs unless $(d, q) = (3, 2)$.
For $(d, q) = (3, 2)$, we have four ${\rm P\Gamma{}L}(U)$-double cosets
in Sym$(\cH)$.
\qed

Similar estimates would follow if one generalized Lemma \ref{lem:detU}
to show that $U$ is uniquely determined in $\bP$ for arbitrary designs
$\bD$ (that is, for $\phi$ that are not permutations). 
The blocks of $\bD$ are then found as 
$\{ \Gamma_\phi(x) \cap U \mid x \in P \setminus U \}$.
In \cite[Corollary 3.2]{Kantor94} it is shown that 
for $d \geq 4$ there are at least $q^{d-2}!$
choices for $\bD$.
Hence, one would obtain the same estimate as in Theorem \ref{thm:low_bnd}
for $d \geq 4$. 


\subsection{Switched symplectic graphs with 4-vertex condition}

We show that in the symplectic case the graphs
$\Gamma_\phi$ satisfy the 4-vertex condition.
Let $\bP$ be ${\rm Sp}_{2d}(q)$, and 
let $V$ be a $2d$-dimensional vector space over $\Ff_q$,
provided with a nondegenerate symplectic form.

{\medskip\footnotesize
The parameters of $\Gamma_0$ are
$v=(q^{2d}-1)/(q-1)$,
$k=q(q^{2d-2}-1)/(q-1)$,
$v-k-1=q^{2d-1}$,
$\lambda=q^2(q^{2d-4}-1)/(q-1)+q-1$,
$\mu=(q^{2d-2}-1)/(q-1)$
and
$\binom{\lambda}{2} - \alpha = \frac12 q^{2d-1}(q^{2d-4}-1)/(q-1)$,
$\beta = \frac12 q(q^{2d-2}-1)(q^{2d-4}-1)/(q-1)^2$,
and those of $\Gamma_\phi$ will turn out to be the same.
\par\smallskip}

\begin{Proposition}
    The graph $\Gamma_\phi$ satisfies the $4$-vertex condition.
\end{Proposition}

\Proof
Let $x,y$ be two vertices of $\Gamma_\phi$.
We show that the number of edges in $\Gamma_\phi(x) \cap \Gamma_\phi(y)$
is independent of $\phi$, and only depends on whether $x,y$ are adjacent
or nonadjacent.
Since $\Gamma_0$ satisfies the 4-vertex condition, $\Gamma_\phi$ does too.

\smallskip 
Count edges $ab$ in $\Gamma_\phi(x) \cap \Gamma_\phi(y)$.
The vertices $x,y,a,b$ are pairwise adjacent, except that
$x$ and $y$ need not be adjacent. We distinguish several cases depending on
which of $x,y,a,b$ are in $U$. Each of the separate counts will be
independent of $\phi$.
If $x \notin U$ then let $X = x^\perp \cap U$.
If $y \notin U$ then let $Y = y^\perp \cap U$.

\paragraph*{Case $x,y,a,b \notin U$.}
In this case adjacencies and counts do not involve $\phi$.

\paragraph*{Case $a,b \in U$.}
Here $a,b$ must be chosen distinct from $x,y$ in case
$x,y \in U$, or distinct from $x$ and in $Y^\phi$ in case
$x \in U$, $y \notin U$ (and the count depends on whether $x \adj y$), or in
$X^\phi \cap Y^\phi$ in case $x,y \notin U$
(and the count depends on whether $X = Y$).
In all cases the count is independent of $\phi$.

\paragraph*{Case $x,y,a \in U$, $b \notin U$.}
For each hyperplane $H$ such that $x,y \in H^\phi$ we count
the $b \in H^\perp \setminus U$ and the $a \in H^\phi$
distinct from $x,y$.

\paragraph*{Case $x,y \in U$, $a,b \notin U$.}
For any two hyperplanes $H,H'$ of $U$ with $x, y \in H^\phi \cap H'^\phi$
count adjacent $a,b$ with $a \in H^\perp \setminus U$
and $b \in H'^\perp \setminus U$.
(The counts will depend on whether $H = H'$, but not on $\phi$.)

\paragraph*{Case $x,a \in U$, $y,b \notin U$.}
For each hyperplane $H$ with $x \in H^\phi$, count
the $a \in H^\phi \cap Y^\phi$ distinct from $x$, and
$b \in H^\perp \setminus U$ adjacent to $y$.
(Here $H = Y$ occurs when $x \adj y$. The counts for $H \ne Y$
do not depend on $H$.)

\paragraph*{Case $x \in U$, $y,a,b \notin U$.}
For any two hyperplanes $H,H'$ with $x \in H^\phi \cap H'^\phi$, count
edges $ab$ with $a \in H^\perp$ and $b \in H'^\perp$
in $y^\perp \setminus (U \cup \{y\})$.
(Here $H = Y$ or $H' = Y$ occur when $x \adj y$.
The counts for $H,H' \ne Y$ do not depend on the hyperplanes chosen
but only on whether $H = Y$ or $H' = Y$ or $H = H'$.)

\bigskip
Finally the least trivial case.

\paragraph*{Case $a \in U$, $x,y,b \notin U$.}
Count $a,H,b$ with $a \in X^\phi \cap Y^\phi$ and $H$
a hyperplane of $U$ on $a$ and $b \in \<x,y,H\>^\perp
\setminus (U \cup \{x,y\})$. The count for $a$ depends
on whether $X = Y$, that for $b$ depends on whether
$H = X$ or $H =Y$ or $H \supseteq X \cap Y$,
but does not otherwise depend on the choice of $H$.

\bigskip
Since all counts were independent of $\phi$, this proves
our proposition.
\qed

\medskip

By Theorem \ref{thm:low_bnd}, this shows that there
are many strongly regular graphs which satisfy the 
4-vertex condition. But we still have to show the 
simplified version of this statement given in 
the introduction as Theorem \ref{thm:4vtxsrgs}.

\medskip

\noindent{\bf Proof of Theorem \ref{thm:4vtxsrgs}.}\quad 
Note that here $v$ refers to a nonnegative integer
as in Theorem \ref{thm:4vtxsrgs} and no longer is the 
number of vertices in $\Gamma_\phi$.

Apply Theorem \ref{thm:low_bnd} for $d=3$ to find at least
$q!$ strongly regular graphs satisfying the 4-vertex condition
on $\tilde{v}$ vertices, for $\tilde{v} = \smash{\frac{q^6-1}{q-1}}$.
Given $v$, there is a prime $q$ between $v^{1/6}$ and $2v^{1/6}$
by Bertrand's postulate. Now $\tilde{v} < 2q^5 < 64 v^{5/6} < v$
for $v > 2^{36}$.
%
%
Checking the prime powers $q$ for $7 \le q \le 64$
one sees that there is a $q$ with $\tilde{v} \le v \le q^6$
for $v \ge 19608$.
One easily verifies the assertion for $v < 19608$ using rank 3 graphs.
%
%
\qed

{\medskip\footnotesize 
Further graphs with the same parameters satisfy the 4-vertex condition.
Additional examples can be obtained by repeated WQH-switching,
see \S\ref{subsec:small} and \cite{Ihringer19}, and there are more examples
among the graphs constructed in \cite{Ihringer17}.
We have not tried (much) to determine precisely which
graphs in \cite{Ihringer17} do satisfy the 4-vertex condition. Similarly,
we do not know when WQH-switching preserves the 4-vertex 
condition.
\par}

%
%

\subsection{Small examples}\label{subsec:small}

\subsubsection*{Examples on 63 vertices}

In \cite{Ihringer20} a large number of strongly regular graphs are found
by applying GM-switching to the ${\rm Sp}_6(2)$ polar graph.
Among these are 280 non-rank-3 strongly \mysqueeze{0.25pt}{regular graphs
with $(v,k,\lambda,\mu) = (63, 30, 13, 15)$ satisfying
the 4-vertex condition. All have $\alpha=30$ and $\beta=45$.
Three of these are among the $\Gamma_\phi$ constructed above.}

We list for each occurring group size the number of examples found.

{\medskip\footnotesize\noindent
\setlength{\tabcolsep}{4pt}
\begin{tabular}{c|cccccccccccccccc}
$|G|$ & 4 & 8 & 16 & 32 & 48 & 64 & 96 & 128 & 192 & 256 & 384 & 512 & 768 & 1344 & 1536 & 4608 \\
\hline
\# & 3 & 16 & 76 & 62 & 1 & 60 & 2 & 30 & 5 & 12 & 3 & 3 & 2 & 1 & 3 & 1
\end{tabular}
\par\medskip}

None of these examples has a transitive group.
We list the orbit lengths in the seven cases with fewer than six orbits.

{\medskip\footnotesize\noindent
\begin{tabular}{c|cccccccccccccccc}
$|G|$ & 768 & 768 & 1344 & 1536 & 1536 (twice) & 4608 \\
\hline
orbits & $3{+}12{+}48$ & $1{+}6{+}24{+}32$ & $7 {+} 56$ & $1{+}6{+}24{+}32$ & $3{+}4{+}8{+}48$ & $3{+}12{+}48$
\end{tabular}
\par\smallskip}

\subsubsection*{Permutations of hyperplanes}

Let $\bP$ be ${\rm Sp}_{2d}(q)$, and let $\phi$ be a permutation
of the set $\cH$ of hyperplanes of $U$.
For $(d, q) = (3, 2)$, $(3, 3)$, $(4, 2)$, the number of double cosets
of ${\rm P\Gamma{}L}(d, q)$ in ${\rm Sym}(\cH)$
is $4$, $252$, and $3374$, respectively, and these are the numbers
of non-isomorphic graphs $\Gamma_\phi$. In each case, exactly one has rank $3$.
None of the others has a transitive group (since $U$ can be recognized).
The pointwise stabiliser of $U$ in ${\rm Aut}(\Gamma_0)$ has size
$N = q^{\binom{d+1}{2}} (q-1)$
and is always contained in ${\rm Aut}(\Gamma_\phi)$. 
Hence, $N$ divides $|{\rm Aut}(\Gamma_\phi)|$.

\bigskip
{\em Case $(d, q) = (3, 3)$.}
Here $N = 1458$.
We list the group sizes for the 251 graphs $\Gamma_\phi$ other than
$\Gamma_0$.

{\medskip\footnotesize\noindent
\setlength{\tabcolsep}{4pt}
\begin{tabular}{c|cccccccccccccc}
$|G|/N$ & 1 & 2 & 3 & 4 & 6 & 8 & 12 & 16 & 18 & 24 & 39 & 54 & 72 & 144 \\
\hline
\# & 172 & 26 & 29 & 6 & 3 & 2 & 2 & 2 & 1 & 1 & 3 & 1 & 2 & 1 
\end{tabular}
\par\medskip}

We list the orbit lengths in the five cases with fewer than six orbits.

{\medskip\footnotesize\noindent
\begin{tabular}{c|ccc}
$|G|/N$ & 39 (thrice) & 72 & 144 \\
\hline
orbits & $13{+}351$ & $1{+}12{+}108{+}243$ & $1{+}12{+}108{+}243$
\end{tabular}
\par}

\bigskip
{\em Case $(d, q) = (4, 2)$.}
Here $N = 1024$.
We list the group sizes for the 3373 graphs $\Gamma_\phi$ other than
$\Gamma_0$.

{\medskip\footnotesize\noindent
\setlength{\tabcolsep}{2.75pt}
\begin{tabular}{c|cccccccccccccccccccc}
$|G|/N$ & 1 & 2 & 3 & 4 & 5 & 6 & 7 & 8 & 12 & 16 & 18 & 21 & 24 & 32 & 56 & 60 & 96 & 192 & 288 & 1344 \\
\hline
\# & 3148 & 85 & 40 & 24 & 4 & 10 & 6 & 26 & 1 & 4 & 1 & 2 & 11 & 2 & 2 & 1 & 2 & 2 & 1 & 1
\end{tabular}
\par\medskip}

We list the orbit lengths in the eight cases with fewer than six orbits.

{\medskip\footnotesize\noindent
\begin{tabular}{c|cccccc}
$|G|/N$ & 12 & 18 & 24 & 56 (twice) \\
\hline
orbits & $3{+}12{+}48{+}192$ & $6 {+} 9 {+} 96 {+} 144$ & $3 {+} 12 {+} 48 {+} 192$ & $1 {+} 14 {+} 112 {+} 128$ 
\end{tabular}

\noindent
\begin{tabular}{c|ccccc}
$|G|/N$ & 60 & 288 & 1344 \\
\hline
orbits & $15 {+} 240$ & $3{+}12{+}48{+}192$ & $7{+}8{+}16{+}224$
\end{tabular}

\par\smallskip}

\subsubsection*{Other polar spaces}

We made the same exhaustive investigation of all permutations $\phi$
for the other choices of $\bP$ in the cases
$(d, q) \in \{ (3, 2), (3, 3), (4, 2) \}$.
The only non-rank-3 examples satisfying the 
$4$-vertex condition occur for ${\rm O}_7(3)$.
Here we obtain $252$ graphs in total, of which one is rank $3$,
and three more satisfy the $4$-vertex condition.
They all have two orbits (of sizes $13{+}351$) and
an automorphism group of size 56862.
All other graphs $\Gamma_\phi$ obtained from $O_7(3)$ have 
more than two orbits.

One might wonder whether
a graph $\Gamma_\phi$ from ${\rm O}_{2d+1}(q)$ satisfies
the $4$-vertex condition if and only if it has 
at most two orbits. And whether a non-rank-3 graph $\Gamma_\phi$
can only satisfy the 4-vertex condition if $\bP$
is ${\rm Sp}_{2d}(q)$ or ${\rm O}_{2d+1}(q)$.

\subsubsection*{Other designs}

There are four 2-$(15, 7, 3)$ designs $\bD$
other than that of the hyperplanes of ${\rm PG}(3,2)$.
We investigated the case where $(d,q) = (4, 2)$ and $\bP$ is ${\rm Sp}_2(8)$,
so that the resulting examples satisfy the 4-vertex condition.
We generated several hundred thousand graphs $\Gamma_\phi$
for each of these designs. 
None of these graphs occurs for two different designs.
We believe our enumeration to be complete.

\medskip

\begin{tabular}{cccc}
 $|{\rm Aut}(\bD)|$ & point orbits & block orbits & \# $\Gamma_\phi$ \\ \hline
 576 & $3{+}12$ & $3{+}12$ & 113519 \\
 168 & $7{+}8$ & $1{+}14$ & 340730 \\
 168 & $1{+}14$ & $7{+}8$ & 328078 \\
 96 & $1{+}6{+}8$ & $1{+}6{+}8$ & 677460
\end{tabular}

\begin{appendix}
\section*{Appendix A --- Details on Ivanov's graphs}
{\footnotesize
In Section \ref{binaryvs} we discussed the graphs
$\Gamma^{(m)}$ from \cite{BrouwerIvanovKlin89}
and $\Sigma^{(m)}$ from \cite{Ivanov94}.
Here we give some more detail on the latter.

For $m \ge 2$, consider $V = \Ffss{2}{2m}$ provided with
the elliptic quadratic form $q(x) = x\subsupr{1}{2} + x\subsupr{2}{2} +
x_1x_2 + x_3x_4 + ... + x_{2m-1}x_{2m}$.
Identify the set of projective points (1-spaces) in $V$
with $V^* = V \setminus \{0\}$.
Let $Q = \{ x \in V^* \mid q(x) = 0 \}$ and let $S$ be the
maximal t.s.~subspace given by $S = \{ x \in V^* \mid
x_1 = x_2 = 0 ~{\rm and}~ x_{2i-1}=0 ~(2 \le i \le m) \}$.
Then $S^\perp = \{ x \in V^* \mid x_{2i-1}=0 \maysplit ~(2 \le i \le m) \}$.
The graph $\Sigma^{(m)}$ has $V$ as vertex set, where two distinct vertices
$v,w$ are adjacent when $v-w \in (Q \cup S^\perp) \setminus S$.
Let $\Tau^{(m)}$ and $\Upsilon^{(m)}$ be the induced
subgraphs on the neighbors (nonneighbors) of the vertex 0.
Put $R = V^* \setminus (Q \cup S^\perp)$.

\medskip\noindent{\bf Proposition.}

(i) For $m \le 4$, the graphs $\Sigma^{(m)}$ are rank $3$,
and are isomorphic to the complement of $VO_{2m}^-(2)$.

(ii) For $m \ge 5$, the automorphism group of $\Tau^{(m)}$ has
two vertex orbits $S^\perp \setminus S$ and $Q \setminus S$,
of sizes $3 \cdot 2^{m-1}$ and $2^{2m-1}-2^m$, respectively.
For $2 \le m \le 4$, the group is rank $3$,
and the graph is the complement of $NO_{2m}^-(2)$.

(iii) For $m \ge 5$, the automorphism group of $\Upsilon^{(m)}$
has two vertex orbits $S$ and $R$ of sizes
$2^{m-1}-1$ and $2^{2m-1}-2^m$, respectively.
For $3 \le m \le 4$, the group is rank $3$,
and the graph is the complement of $O_{2m}^-(2)$.

(iv) The $\lambda$- and $\mu$-graphs in $\Upsilon^{(m)}$
and the $\mu$-graphs in $\Tau^{(m)}$ are all regular
of valency $2^{m-2}(2^{m-2}+1)$.
In particular, $\Upsilon^{(m)}$ satisfies the 4-vertex condition.

(v) The $\lambda$-graphs in $\Tau^{(m)}$ have vertices of valencies
in $0$, $2^{2m-4}-2^m$, $2^{2m-4}$, $2^{2m-3}-2^m$.
Edges not in a line contained in $Q$ have $\lambda$-graphs
with a single isolated vertex and $\lambda-1$ vertices of valency $2^{2m-4}$.
For edges in a line contained in $Q$ the $\lambda$-graphs have
a single vertex with valency $2^{2m-3}-2^m$, and
$2^{m-3}-1$ vertices with valency $2^{2m-4}-2^m$,
and the remaining $2^{2m-3}+2^{m-3}$ vertices have valency $2^{2m-4}$.
In particular, $\Tau^{(m)}$ satisfies the 4-vertex condition,
with $\alpha = 2^{2m-5}(2^{2m-3}+2^{m-2}-1)$ and
$\beta = \frac12 \mu \mu' = 2^{2m-4} (2^{m-2}+1)^2$.

(vi) The local graph of $\Upsilon^{(m)}$ at a vertex $s \in S$
is isomorphic to $\Sigma^{(m-1)}$.

\medskip\noindent\Proof
(i)--(iii) This is clear, and can also be found in \cite{Ivanov94}.

(iv)-(v) (the part about $\Tau^{(m)}$):

Let $(v,w) = q(v+w)-q(v)-q(w)$ be the symmetric bilinear form
belonging~to~$q$.
Let $X = (Q \cup S^\perp) \setminus S$. Then $\Tau^{(m)}$
is the graph with vertex set $X$, where two vertices $x,y$ are adjacent
when the projective line $\{x,y,x+y\}$ they span is contained in $X$.
If at least one of $x,y$ is in $S^\perp \setminus S$, then this is
equivalent to $(x,y)=1$. If both are in $Q \setminus S$, then this is
equivalent to ($(x,y)=0$ and $x+y \notin S$) or ($(x,y)=1$ and
$x+y \in S^\perp \setminus S$).

Let $x,y,z$ be pairwise adjacent vertices. The valency $c$ of $z$ in
the $\lambda$-graph $\lambda(x,y)$ is the number of common neighbors
of $x,y,z$. Distinguish several cases.

If $z = x+y$, then if $x,y,z \in Q$ we find
$c = | \{x,y\}^\perp \cap (Q \setminus S) | - 3 = 2^{2m-3}-2^m$.
If $z = x+y$ and at least one of $x,y,z$ lies in $S^\perp$, then $c = 0$.

Now let $z \ne x+y$.
The claims are true for $m \le 4$. Let $m \ge 5$ and use induction on $m$.
Choose coordinates so that $x,y,z$ have final coordinates $00$
and let $x',y',z'$ be these points without the final two coordinates.
If they have $c'$ common neighbors $w'$ in $\Tau^{(m-1)}$,
then we find $2c'$ common neighbors $w = (w',0,*)$.
Moreover (since $x,y,z$ are linearly independent), we find
$2^{2m-5}$ common neighbors $(w',1,q'(w'))$ in $Q$, where $w'$ runs through
all vectors with the desired inner products with $x',y',z'$.
Altogether $c = 2c'+2^{2m-5}$, as claimed.

For the $\mu$-graphs the argument is similar and simpler:
by the definition of adjacency three dependent vertices are pairwise adjacent,
so that the case $z = x+y$ does not occur here.

(iv) (the part about $\Upsilon^{(m)}$):
Let $Y = V^* \setminus X$. Then $\Upsilon^{(m)}$ is the graph
with vertex set $Y$, where two vertices $x,y$ are adjacent when
the projective line $\{x,y,x+y\}$ they span is not contained in $Y$.
The same argument as before yields the valencies of the
$\lambda$- and $\mu$-graphs.

(vi) Consider the graph $\Sigma^{(m)}$.
The nonneighbors $z$ of 0 that are neighbors of $s$ are the vertices
of the form $z=s+b$ with $z \in S \cup R$
and $b \in (Q \cup S^\perp) \setminus S$.
It follows that $s+z \in Q \setminus s^\perp$.
Let $s = (0\ldots 01)$, then $Q \setminus s^\perp$ can be identified with
$W = \Ffss{2}{2m-2}$ via $w \to i(w)=(w,1,\bar{q}(w))$
for $w \in \Ffss{2}{2m-2}$ and $\bar{q}(w)$ determined by $q(i(w))=0$.
The local graph of $\Upsilon$ at $s$ can be identified with
the graph with vertices $w$, where $w,w'$ are adjacent when the line joining
$i(w),i(w')$ has third point $(w+w',0,*) \in (Q \cup S^\perp) \setminus S$,
that is, the line joining $w,w'$ has third point $w''=w+w'$ satisfying
$w'' \notin T$ and $(\bar{q}(w'')=0$ or $w'' \in T^\perp)$ where
$T = \{ w \in W \mid w_1=w_2=w_3=w_5=...=w_{2m-3}=0 \}$.
But this is $\Sigma^{(m-1)}$.
\qed

\par}

\end{appendix}

\bigskip
\paragraph*{Acknowledgment} The second author is supported by a 
postdoctoral fellowship of the Research Foundation -- Flanders (FWO).


\begin{thebibliography}{99}

\bibitem{BannaiHaoSong90}
E. Bannai, S. Hao \& S.-Y. Song,
{\it Character tables of the association schemes of finite
orthogonal groups acting on the nonisotropic points},
J. Comb. Th. (A) {\bf 54} (1990) 164--200.

\bibitem{Brouwer16}
A. E. Brouwer,
{\it Strongly regular graphs from hyperovals}, \quad\verb|https://|
\verb|www.win.tue.nl/~aeb/preprints/hhl.pdf|,
accessed on 2021-02-21.

\bibitem{BCN}
A. E. Brouwer, A. M. Cohen \& A. Neumaier,
{\it Distance-regular graphs}, Springer, Heidelberg, 1989.

\bibitem{BrouwerIvanovKlin89}
A. E. Brouwer, A. V. Ivanov \& M. H. Klin,
{\it Some new strongly regular graphs},
Combinatorica {\bf 9} (1989) 339--344.

\bibitem{BrouwerVanLint84}
A. E. Brouwer \& J. H. van Lint,
{\it Strongly regular graphs and partial geometries},
pp. 85--122 in: Enumeration and design (Waterloo, Ont., 1982),
Academic Press, 1984.


\bibitem{BrouwerVM21}
A. E. Brouwer \& H. Van Maldeghem,
{\it Strongly regular graphs},
Cambridge Univ. Press, Cambridge, 2022.

\bibitem{Cameron74}
P. J. Cameron,
{\it Partial quadrangles},
Quart. J. Math. Oxford, {\bf 25(3)} (1974), 1--13.

\bibitem{CameronGoethalsSeidel78}
P. J. Cameron, J. M. Goethals \& J. J. Seidel,
{\it Strongly regular graphs having strongly regular subconstituents},
J. Algebra {\bf 55} (1978) 257--280.



\bibitem{Debroey78}
I. Debroey,
{\it Semi partiële meetkunden},
Ph.~D.~thesis, University of Ghent, 1978.

\bibitem{DebroeyThas78}
I. Debroey \& J. A. Thas,
{\it On semipartial geometries},
J. Comb. Th. (A) {\bf 25} (1978) 242--250.

\bibitem{Delsarte72a}
Ph. Delsarte, 
{\it Weights of linear codes and strongly regular normed spaces}, 
Discr. Math. {\bf 3} (1972) 47--64.

\bibitem{DempwolffKantor08}
U. Dempwolff \& W. M. Kantor,
{\it Distorting symmetric designs},
Des. Codes Cryptogr. {\bf 48} (2008) 307--322.




\bibitem{HestenesHigman71}
M. D. Hestenes \& D. G. Higman,
{\it Rank 3 groups and strongly regular graphs},
pp.~141--159 in: Computers in algebra and number theory
(Proc. New York Symp., 1970), G. Birkhoff \& M. Hall jr (eds.),
SIAM-AMS Proc., Vol IV, Providence, R.I., 1971.

\bibitem{Higman71}
D. G. Higman,
{\it Partial geometries, generalized quadrangles and strongly regular graphs},
pp. 263--293 in:
Atti del Convegno di Geometria Combi\-natoria e sue Applicazioni
(Univ. Perugia, Perugia, 1970), Ist. Mat., Univ. Perugia, Perugia (1971).

\bibitem{Hill73}
R. Hill,
{\it Caps and groups},
pp. 389--394 in:
Colloquio Internazionale sulle Teorie Combinatorie (Rome, 1973),
Tomo II, Atti dei Convegni Lincei, No. 17, Accad. Naz. Lincei, Rome, 1976.

\bibitem{HirschfeldThas80}
J. W. P. Hirschfeld \& J. A. Thas,
{\it Sets of type $(1,n,q+1)$ in $PG(d,q)$},
Proc. London Math. Soc. (3) {\bf 41} (1980) 254--278.

\bibitem{HuangHuangLin09}
T. Huang, L. Huang \& M.-I. Lin, 
{\it On a class of strongly regular designs and quasi-semisymmetric designs},
pp. 129–153 in: Recent developments in algebra and related areas,
Proceedings Conf. Beijing 2007, Chongying Dong et al. (eds.), Adv. Lect. Math.
(ALM) 8, Higher Education Press and Int. Press, Beijing-Boston, 2009.

\bibitem{Ihringer17}
F. Ihringer,
{\it A switching for all strongly regular collinearity graphs from
polar spaces},
J. Algebr. Comb. {\bf 46} (2017), 263--274.

\bibitem{Ihringer19}
F. Ihringer \& A. Munemasa,
{\it New strongly regular graphs from finite geometries via switching},
Linear Algebra Appl. {\bf 580} (2019), 464--474.

\bibitem{Ihringer20}
F. Ihringer,
{\it Switching for Small Strongly Regular Graphs},\quad\quad
{\tt arXiv:\allowbreak 2012.08390v1} (2020). 

\bibitem{Ivanov89}
A. V. Ivanov,
{\it Non rank 3 strongly regular graphs with the 5-vertex condition},
Combinatorica {\bf 9} (1989) 255--260.
%
%

\bibitem{Ivanov94}
A. V. Ivanov,
{\it Two families of strongly regular graphs with the 4-vertex condition},
Discr. Math. {\bf 127} (1994) 221--242.

\bibitem{Kantor86}
W. M. Kantor,
{\it Some generalized quadrangles with parameters $(q^2, q)$},
Math. Z. 192 (1986) 45--50.

\bibitem{Kantor94}
W. M. Kantor,
{\it Automorphisms and isomorphisms of symmetric and affine designs},
J. Alg. Comb. {\bf 3} (1994) 307--338.

\bibitem{Kaski-et-al12}
P. Kaski, M. Khatirinejad \& P. R. J. Östergård,
{\it Steiner triple systems satisfying the 4-vertex condition},
Des. Codes Cryptogr. {\bf 62} (2012) 323--330.

\bibitem{Klin-et-al05}
M. Klin, M. Meszka, S. Reichard \& A. Rosa,
{\it The smallest non-rank 3 strongly regular graphs
which satisfy the 4-vertex condition},
Bayreuther Mathematische Schriften {\bf 74} (2005) 145--205.

\bibitem{KlinPech}
M. Klin \& C. Pech, May 2008, unpublished notes.


\bibitem{PayneThas84}
S. E. Payne \& J. A. Thas,
{\it Finite generalized quadrangles},
Research Notes in Mathematics, 110. 
Pitman (Advanced Publishing Program), Boston, MA, 1984. vi+312 pp. 

\bibitem{PechPech19}
C. Pech \& M. Pech,
{\it On a family of highly regular graphs by Brouwer, Ivanov, and Klin},
Discr. Math. {\bf 342} (2019) 1361--1377.

\bibitem{Reichard00}
S. Reichard,
{\it A criterion for the $t$-vertex condition on graphs},
J. Comb. Th. (A) {\bf 90} (2000) 304--314.

\bibitem{Reichard15}
S. Reichard,
{\it Strongly regular graphs with the 7-vertex condition},
J. Algebr. Comb. {\bf 41} (2015) 817--842.

\bibitem{Sims-unpub}
C. C. Sims,
{\it On graphs with rank 3 automorphism groups},
unpublished, 1968.


\bibitem{WangQiuHu19}
W. Wang, L. Qiu \& Y. Hu,
{\it Cospectral graphs, GM-switching and
regular rational orthogonal matrices of level $p$},
Lin. Alg. Appl. {\bf 563} (2019) 154--177.

\bibitem{Wilbrink}
H. A. Wilbrink, unpublished, 1982.


\end{thebibliography}
\end{document}